\documentclass{amsart}
\usepackage{amssymb,latexsym}
\usepackage{amsmath}
\usepackage{graphicx}
\usepackage{color}

\newtheorem{thm}{Theorem}[section]
\newtheorem{cor}[thm]{Corollary}

\newtheorem{lemma}[thm]{Lemma}

\newcommand{\bdd}{\mbox{$\partial$}}
\newcommand{\aaa}{\mbox{$\alpha$}}
\newcommand{\bbb}{\mbox{$\beta$}}

\begin{document}



\thanks{Research of the first author is partially supported by an NSF grant.}

\title{Uniqueness of bridge surfaces for 2-bridge knots}

\author{Martin Scharlemann}
\address{\hskip-\parindent
Martin Schrlemann\\
Mathematics Department \\
University of California, Santa Barbara \\
Santa Barbara, CA 93117, USA}
\email{mgscharl@math.ucsb.edu}

\author{Maggy Tomova}
\address{\hskip-\parindent
Maggy Tomova\\
Mathematics Department \\
University of Iowa \\
 Iowa city, IA 52240, USA}
\email{mtomova@math.uiowa.edu}

\date{\today}

\begin{abstract}
 
Any $2$-bridge knot in $S^3$ has a bridge sphere from which any other bridge surface can be obtained by stabilization, meridional stabilization, perturbation and proper isotopy.  

\end{abstract}

\maketitle

\section{Introduction}

Establishing the uniqueness of Heegaard splittings for certain 3-manifolds has been an 
interesting and surprisingly difficult problem. One of the earliest 
known results was that of Waldhausen \cite{wald68} who proved that $S^3$ has a 
unique Heegaard splitting up to stabilization. In \cite{BO83}, Bonahon 
and Otal proved that the same is true of lens splaces (manifolds 
with a genus one Heegaard surface). A later proof \cite{RS1} made use of the fact that any two weakly incompressible Heegaard  splittings of a manifold can 
be isotoped to intersect in a nonempty collection of curves that are  essential on both Heegaard surfaces. 

There is an analogue to Heegaard splitting in the theory of links in $3$-manifolds. (By link, we include the possibility that $K$ has one component, i. e. a knot is a link.)   Consider a link $K$  in a closed orientable $3$-manifold $M$ with a Heegaard surface $P$ (i.e. $M=A \cup _P B$ where $A$ and $B$ are handlebodies) and require that each arc of $K-P$ is $P$-parallel in $M - P$.  We say that $K$ is in {\em bridge position} with respect to $P$ and 
that $P$ is a {\em bridge surface} for the pair $(M,K)$.  Beyond the philosophical analogy between Heegaard splittings for $3$-manifolds and bridge surfaces for links in $3$-manifolds, notice that there is also this precise connection:  If $P$ is a bridge surface for a link $K$ in $M$, then the cover $\hat{P}$ of $P$ in the $2$-fold branched cover $\hat{M}$ of $M$ is a Heegaard surface for the manifold $\hat{M}$.  

Questions about the structure of Heegaard splittings on $3$-manifolds often have analogies with questions about bridge surfaces.  For example, it is natural to ask whether there are pairs $(M,K)$ that have a unique bridge surface, up to some obvious geometric operations analogous to Heegaard stabilization. In \cite {HS2} Hayashi and Shimokawa proved that this is true for bridge surfaces of the unknot. We will show that the same is true for bridge surfaces of $2$-bridge knots. (And presumably for $2$-bridge links as well, though we do not pursue that here, because of the technical obstacle that the theory in \cite{STo} so far has not been explicitly extended to $3$-manifolds with non-empty boundary. Compare \cite{RS2} to \cite{RS1}.)  This result can viewed as the analogue for bridge surfaces of the result of Bonahon and Otal mentioned above.  Our approach will be analogous to that of \cite{RS1}, working from the central result of \cite{STo}: in the absence of incompressible Conway spheres, two c-weakly incompressible bridge surfaces can be properly isotoped to intersect in a non-empty collection of closed curves, each of which is essential (including non-meridional) in both surfaces.  

\section{Definitions and Notation}

If $X$ is any subset of a 3-manifold $M$ and $K$ is a 1-manifold 
properly embedded in $M$, let $X_K=X-K$. A disk $D\subset M$ that
meets $K$ exactly once is called a {\em punctured disk}. If $F$ is
an embedded surface in $M$ transverse to $K$, a simple 
closed curve on $F_K$ is
{\em essential} if it doesn't bound a disk or a punctured disk on $F_K$. An 
embedded
disk $D \subset M_K$ is a {\em compressing disk} for $F_K$ if $D \cap 
F_K
=\bdd D$ and $\bdd D$ is an essential curve in $F_K$. A {\em cut-disk}
for $F_K$ is a punctured disk $D^c$ in $M_K$ such that $D^c \cap 
F_K =\bdd D^c$ and $\bdd D^c$ is an essential curve in $F_K$. A 
possibly punctured disk $D^*$ that is either a cut disk or a compressing disk will be
called a {\em c-disk} for $F_K$. The surface $F_K$ is called essential if it has no 
compressing disks (it may have cut-disks), it is not a sphere that 
bounds a ball in $M_K$ and it is not $\bdd$-parallel in $M-\eta(K)$ where $\eta(K)$ is a 
regular open neighborhood of $K$.

A properly embedded arc $\alpha \subset F_K$ is {\em inessential} if there 
is a disk on $F_K$ whose boundary is the endpoint union of $\alpha$ 
and a subarc of $\bdd F$. Otherwise $\alpha$ is {\em essential}. A 
$\bdd$-{\em compressing disk} for $F_K$ is an embedded disk $D \subset M$ 
with an interior disjoint from $F_K$
such that $\bdd D$ is the endpoint union of an essential arc of $F_K$ 
and an arc lying in $\bdd M$.

Any term describing the compressibility of a surface can be extended
to account not only for compressing disks but also c-disks. A surface 
in $M$ that is transverse to $K$ will be called {\em c-incompressible} 
if it has no c-disks. A surface $F$ in $M$ is called a {\em splitting 
surface} if $M$ can be writen as the union of two 3-manifolds along 
$F$. If 
$F$ is a
splitting surface for $M$, we will call $F_K$
{\em c-weakly incompressible}
if any pair of c-disks for $F_K$ on opposite sides of the surface
intersect. If $F_K$ is not c-weakly incompressible, it is {\em 
c-strongly
compressible}. 

A properly embedded collection of arcs $T = \cup_{i=1}^n \alpha_i$ in a compact $3$-manifold is called {\em boundary parallel} if there is a collection $E = \cup_{i=1}^n E_i$ of embedded disks, so that, for each $1 \leq i \leq n$, $\bdd E_i$ is the end-point union of $\alpha_i$ and an arc in the boundary of the $3$-manifold.  A standard cut-and-paste arguments shows that if there is such a collection, there is one in which all the disks are disjoint.   If the manifold is a handlebody $A$, the arcs are called {\em bridges} and disks of parallelism are called {\em bridge disks}. 
Let $M$ be a closed irreducible 3-manifold and let $P$ be a Heegaard 
surface for $M$ decomposing the manifold into handlebodies $A$ and $B$. 
A link $K$ is in bridge position with respect to 
$P$ if each collection of arcs $A \cap K$ and $B\cap K$ is a collection of bridges.
We say that $P$ is a bridge surface for 
the pair $(M,K)$ and the triple $(M; P, K)$ is a bridge presentation of $K \subset M$.

Two disjoint surfaces $F, S \subset M$ transverse to $K$ will be called parallel if they 
cobound a product region and all arcs of the link in that region can 
be isotoped to be vertical with respect to the product structure. $F$ 
is properly isotopic to $S$ if there 
is an isotopy from $F$ to $S$ so that $F$ remains transverse to $K$ throughout the isotopy, 
i.e the isotopy of $F_K$ to $S_K$ is proper in $M_K$.  Unless otherwise stated, all isotopies will be proper isotopies.

\section{New bridge surfaces from old}

Given a bridge surface $P$ for $(M,K)$, 
it is easy to construct more complex bridge surfaces for $(M,K)$ from 
$P$. There are three straightforward ways to do this. The first is easiest: simply add a trivial 1-handle to one of the handlebodies, say $A$. This creates 
a dual 1-handle in $B$. The new bridge surface, 
$P'$ is said to be {\em stabilized} and it is characterized by the presence 
of compressing disks for $P'$, one in $A$ and one in $B$, that intersect in exactly one point.

\medskip

A second way to construct a more complicated bridge surface is almost as easy to see:  Suppose there are a pair of bridge disks $E_A \subset A$ and $E_B \subset 
B$ so that the arcs $E_A \cap P$ and $E_B  \cap P$ intersect precisely at one end.  Then $K$ is said to be {\em perturbed} with respect to $P$ (and vice versa), and $E_A, E_B$ are called {\em cancelling disks} for $K$.   (This is one of two cases of the notion of ``cancellable" bridges, as defined by Hayashi and Shimokawa in \cite{HS2}.  The other case occurs when a component of $K$  is in 1-bridge position, and both bridges, and so a whole component, can be simultaneously isotoped into the bridge surface.)  The word perturbed is used because one way a bridge presentation with this property can be obtained is by starting with any bridge presentation for $K$ and perturbing $K$ near a point of $K \cap P$, introducing a minimum and an adjacent maximum.  The following lemma shows this is in some sense the only way in which a perturbed link can arise.

\begin{lemma}  Suppose $K$ is perturbed with respect to the bridge surface $P$.  Then there is a knot $K'$ in bridge position with respect to $P$, such that $|K' \cap P| = |K \cap P| - 2$ and $K$ is  properly isotopic to the knot obtained from $K'$ by  introducing a minimum and an adjacent maximum near a point of $K' \cap P$.  
\end{lemma}

\begin{proof}  Let $E_A, E_B$ be the cancelling bridge disks, intersecting $P$ in arcs $\aaa$ and $\bbb$ respectively, so that $\aaa \cap \bbb = E_A \cap E_B$ is a single point $p \in P$, an end point of both $\aaa$ and $\bbb$.  A standard cut-and-paste argument shows that there is a disjoint collection of bridge disks for $K \cap A$ so that the collection contains $E_A$. In fact

{\bf Claim:}  There is a disjoint collection $\Delta_A$ of bridge disks for $K \cap A$ so that $E_A \in \Delta_A$ and $\Delta_A \cap \bbb = \bdd \bbb$.  

We begin with a disjoint collection and redefine it so as to eliminate all intersection points with the interior of $\bbb$.  The proof is by induction on the number of points in $\Delta_A \cap interior(\bbb)$.  If the intersection is empty, there is nothing to prove.  Otherwise, suppose that $q$ is the closest point of $\Delta_A \cap \bbb$ to $p$ in $interior(\bbb)$, and let $\bbb'$ be the subsegment of $\bbb$ between $q$ and $p$.  Suppose $E' \neq E_A$ is the bridge disk containing $q$.  Then a regular neighborhood of $E' \cup \bbb' \cup E_A$ has boundary consisting of two disks - one parallel to $E'$ and the other a new bridge disk for the bridge $E' \cap K$ that is disjoint from all other bridge disks and intersects $\bbb$ in one fewer point.  This provides the inductive step, establishing the claim.

Following the claim, let $E' \neq E_A$ be the bridge disk in $\Delta_A$ that is incident to the opposite end of $\bbb$ from $p$; following the claim $E'$, like $E_A$, is disjoint from the interior of $\bbb$. Use $E_B$ to (non-properly) isotope the arc $K \cap E_B$ to $\bbb$ and push it through $P$.  This reduces the number of points in $K \cap P$ by two, but $P$ is still a bridge surface for the knot.  It's clear that $K \cap B$ still consists of bridges, since all we've done is remove one.  The change in $K \cap A$ is to attach the  bridge disk $E'$ to $E_A$ by a band, and the result is clearly still a disk.  It's easy to see that the original positioning of $K$ is properly isotopic to a perturbation of the new positioning of $K$ with respect to $P$.
\end{proof}  

\medskip

Here is a third way to produce a new bridge surface for $(M,K)$, called {\em meridional stabilization}.  Begin with a bridge presentation $M=A \cup_P B$ of $K$ and suppose there is a component $K_0$ of $K$ that is not in $1$-bridge position with respect to $P$.  Let $\bbb$ be a bridge in $K_0 \cap B$ and let $A'$ be the 
union of $A$ together with a neighborhood of $\bbb$.  Let $P' = \bdd A'$ and let $B'$ be the closed complement of $A'$ in $M$.  The decomposition $M=A' \cup_{P'} B'$ is a Heegaard splitting, indeed a stabilization of $M=A \cup_P B$ since a meridian for $A'$ dual to $\bbb$ intersects the remnants of a bridge disk for $\bbb$ in $B'$ in a single point.  Moreover, $K$ is in bridge position with respect to $P'$.  It is obvious that $K \cap B'$ is a collection of bridges, since $K \cap B$ was.  And the new component of $K \cap A'$ has, as a bridge disk, the union of two bridge disks of $K \cap A$ attached together by a band running along $\bbb$.  

\begin{lemma}  A bridge surface $P'$ for $K$ is meridionally stabilized if and only if there is a cut-disk in $A'$ and a compressing disk in $B'$ (or vice versa) that intersect in exactly one point.
\end{lemma}

\begin{proof}  If $P'$ is constructed by meridional stabilization, as described above, then, as we have seen, a meridian disk in $A'$ dual to $\bbb$ is a cut disk for $A'$ that intersects the remnants of a bridge disk for $\bbb$ in a single point. 

Conversely, suppose there is a cut disk $E_A \subset A'$ for $A'$ and a compressing disk $E_B \subset B'$ that intersect in a single point.  Then $P'$ is the stabilization of the Heegaard surface $P$ obtained by cutting $A'$ along $E_A$.  

{\bf Claim:}  $K$ is in bridge position with respect to $P$.

A standard cut and paste argument shows that the bridge disks for $K \cap B'$ can be taken to be disjoint from $E_B$.  They can also be taken to be disjoint from $\bdd E_A$, for any time a bridge disk for $K$ in $B'$ crosses $\bdd E_A$, one can reroute it around $\bdd E_B$, adding a copy of the disk $E_B$ to the bridge disk, to get a bridge disk which intersects $\bdd E_A$ fewer times (see Figure \ref{fig:meridstab}).  Once all bridge disks for $K \cap B'$ are disjoint from $E_A$, they persist when $P'$ is surgered along $E_A$.  So all components of $K \cap B$ have bridge disks, except possibly the new bridge $\bbb$ that is produced in $B$, the bit of $K$ that runs from one copy of $E_A$ (after the cut) to the other.  But $E_B$ itself provides a bridge disk for $\bbb$.  

 \begin{figure}[tbh]
    \centering
    \includegraphics[scale=0.5]{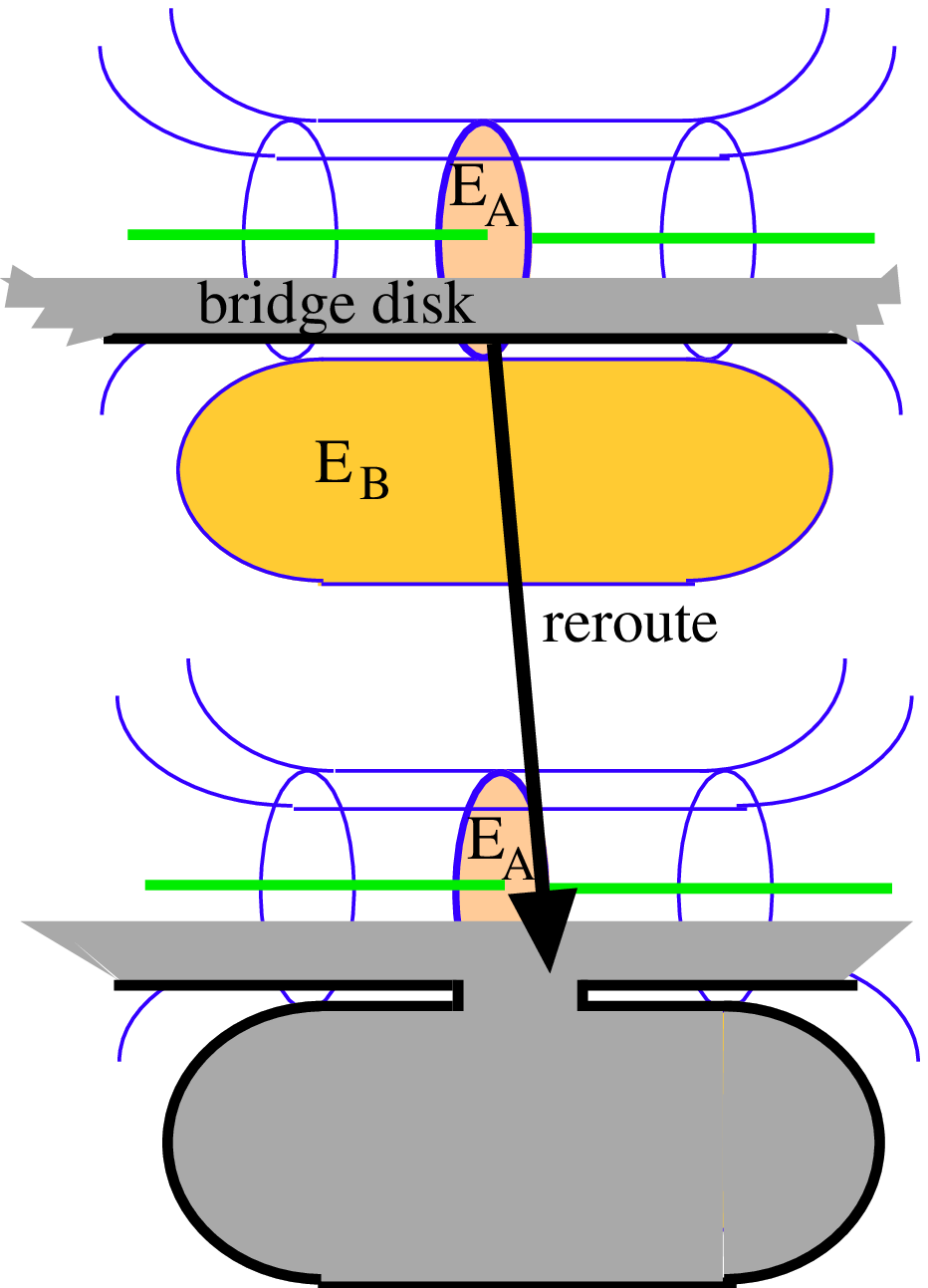}
    \caption{} \label{fig:meridstab}
    \end{figure}

A similar argument exhibits bridge disks in $A$:  A standard cut and paste argument shows that there is a complete collection of bridge disks for $K \cap A'$ that intersects $E_A$ in a single arc, running from the point $K \cap E_A$ to $\bdd E_A$.  When $A'$ is cut apart by $E_A$ to produce $A$, the bridge disk for the component of $K \cap A'$ that intersects $E_A$ is divided by this arc into bridge disks for the two resulting components of $K \cap A$, establishing the claim.

With the claim established, it is easy to see that $P'$ is a meridional stabilization of $P$ along 
$\bbb$.  
\end{proof}

Here is yet a fourth way to construct one bridge surface from another.  It will be useful here to extend, in an obvious way, the definition of bridge surface to links in compact orientable $3$-manifolds with boundary.  Suppose $M$ is a compact orientable $3$-manifold.  A connected closed surface $P \subset M$ is a bridge surface for $K \subset M$ if $P$ is a Heegaard surface for $M$ (that is, the complement of $P$ consists of two compression bodies $C_1, C_2$ and $P = \bdd_+ C_i, i = 1, 2$) and $K$ intersects each complementary compression body in a collection of boundary parallel arcs.  

With that clarifying extension, suppose $K_-$ is a link (possibly empty) in a $3$-manifold $N$ that has a torus boundary component $\bdd_0 N$.  Let $P$ be a bridge surface for $K_-$ in $N$; that is, $P$ divides $N$ into two compression bodies, and $K_-$ intersects each of them in a collection of boundary-parallel arcs.  Fill $\bdd_0 N$ with a solid torus $W$ whose core is a new curve $K_0$.  Then $P$ still divides $M =  N \cup_{\bdd_0 N} W$ into two compression bodies and $K_-$ still intersects each compression body in a collection of boundary-parallel arcs.  Moreover, the core curve $K_0$ is isotopic in $W$ to a curve on $\bdd W = \bdd_0 N$, so $K_0$ is isotopic in $M$ rel $K_-$ to a curve on $P$.  Perturbing $K_0$ slightly makes $P$ a bridge surface for all of $K = K_- \cup K_0$ in $M$.  If a component of a link $K$ in bridge position with respect to $P$ in $M$ can be constructed in this way, then we say that the component is {\em removable}.

\begin{lemma}  \label{lemma:removable}  Suppose $P$ is a bridge surface for a link $K \subset M$.  Then a component $K_0$ of $K$ is removable if and only if $K_0$ can be isotoped rel $K_- = K - K_0$ so that $K_0$ lies on $P$ and there is a meridian disk of one of the two compression bodies that is disjoint from $K_-$ and intersects $K_0 \subset P$ in a single point.
\end{lemma}

\begin{proof}  One direction is fairly straightforward:  if $K_0$ is removable then, in the construction above, $K_0$ can be isotoped to a longitude of $\bdd W$, ie to a curve in $\bdd W$ that intersects a meridian disk $\mu$ of $W$ in a single point.  That is, the wedge of circles $K_0 \vee \bdd \mu \subset \bdd W = \bdd_0 N$.  Let $C$ be the compression body of $N - P$ on which $\bdd_0 N  = \bdd W$ lies.  Then, using the structure of the compression body, there is a proper embedding of $(K_0 \vee \bdd \mu) \times I$ into $C - K_-$, with one end of $(K_0 \vee \bdd \mu) \times I$ on $\bdd W$ and the other end on $P$.  The end on $P$ then describes an embedding of $K_0$ into $P$ that intersects the meridian disk $\mu \cup (\bdd \mu \times I)$ of the compression body $C \cup_{\bdd_0 N} W$ in a single point.  

The other direction uses the ``vacuum cleaner trick":  Suppose that $P$ is a bridge surface for a link $K$ in $M$, that a component $K_0$ of $K$ has been isotoped rel $K_-$ to lie on $P$, and that $\mu$ is a meridian disk for one of the complementary compression bodies $C$ so that $\mu$ is disjoint from $K_-$ and $\mu$ intersects $K_0$ is a single point.  Picture the dual $1$-handle to $\mu$ in $C$ as a vacuum-cleaner  hose, and use it to sweep up all of $K_0 - \eta(\bdd \mu) \subset P$.  Afterwards, $\mu$ is the meridian of a solid torus that is a boundary-summand of $C$, a solid torus for which $K_0$ is a longitude.  Push $K_0$ to the core of this solid torus and remove a thin tubular neighborhood $W$ of $K_0$ from the solid torus.  This changes the solid torus to $torus \times I$, with the result that $C_- = C - W$ is still a compression body.  Moreover, $K_- \cap C_-$ remains a collection of boundary-parallel arcs.  
\end{proof}

\begin{lemma} \label{lemma:remstab}   If a bridge surface for $K$ is stabilized then any $1$-bridge component of $K$ is removable.  

Somewhat conversely, suppose a component $K_0$ for $K$ is removable, with $P$, $K$, $K_0$ and meridian disk $\mu$ as defined in the proof of Lemma \ref{lemma:removable} above.  Suppose further that there is a meridian disk $\lambda$ for the other compression body so that $\lambda$ is disjoint from $K_-$ and $|\mu \cap \lambda| = 1$.  Then $P$ is stabilized.
\end{lemma}

\begin{proof}  Suppose a bridge surface $P$ for $K$ is obtained by stabilizing the bridge surface $P'$ for $K$, and suppose $K_0$ is a $1$-bridge component of $K$. Let $C_1, C_2$ be the compression body complementary components of $P'$.  That is, $|P \cap K_0| = |P' \cap K_0| = 2$, and $P'$ divides $K$ into two boundary-parallel arcs $\tau_i = C_i \cap K, i = 1, 2$.  Let $D_1, D_2$ be bridge disks for $\tau_1, \tau_2$ in $C_1, C_2$ respectively.  By general position, we can assume that the arcs $D_1 \cap P, D_2 \cap P$ have interiors that are disjoint near their end points (though there may be many intersections of their interiors away from the end points).  Stabilize $P'$ to $P$ by attaching a $1$-handle to $C_2$ via an arc $\alpha$ in $D_1$ near and parallel to $\tau_1 \subset \bdd D_1$. Then $D_2$ together with the rectangle in $D_1$ lying between $\alpha$ and $\tau_1$ describes an isotopy of $K_0$ to $P'$.  A cocore of the $1$-handle that was attached to $C_2$ is a meridian for one of the stabilized compression bodies.  Via Lemma \ref{lemma:removable}, $\mu$ exhibits $K_0$ as a removable component of $K$ for the splitting surface $P$.

Now consider the other direction, with meridian disks $\mu \subset C_1$,  $\lambda \subset C_2$, component $K_0 \subset P$ and $|K_0 \cap \mu| = 1 = |\lambda \cap \mu|$ as given in the statement of the lemma.  By general position, we can assume that $K_0$ and $\lambda$ do not intersect near $\mu$.  Move $K_0$ into $1$-bridge position by pushing a small segment of $K_0$ into the interior of $C_2$ near $\mu$ and the interior of the rest of $K_0$ into the interior of $C_1$.  Then $K_0$, hence all of $K$, is disjoint from both meridian disks $\lambda$ and $\mu$, which then exhibit that $P$ is stabilized.
\end{proof}

\begin{figure}[tbh]
    \centering
    \includegraphics[scale=0.6]{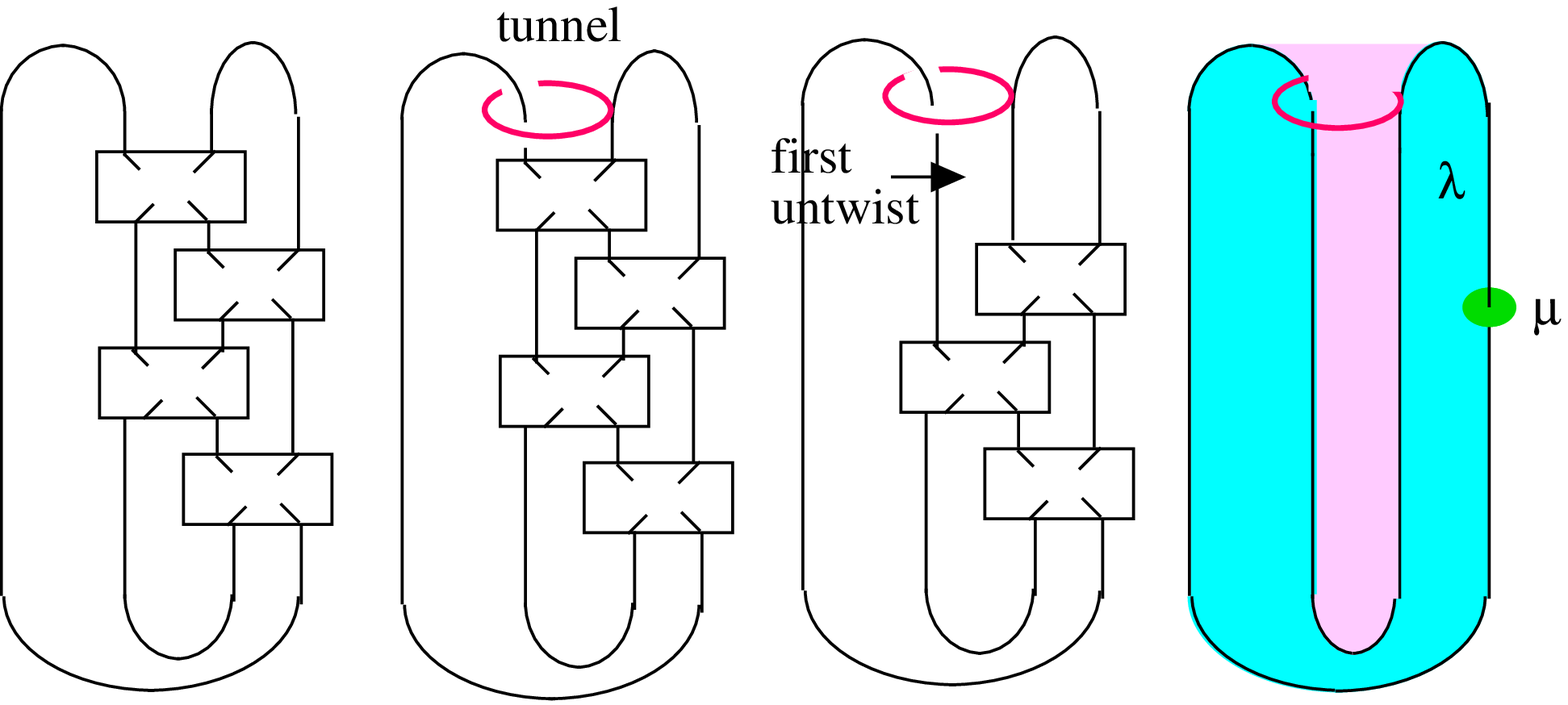}
    \caption{} \label{fig:twobridge}
    \end{figure}

{\bf Example:}  Suppose $K$ is a $2$-bridge knot in $S^3$ and $P$ is a Heegaard surface for the complementary $3$-manifold $N = S^3 - \eta(K)$.  Then either $P$ is stabilized or it is the boundary of a regular neighborhood of the union of the knot and a single arc, and the arc is one of six standard types (see \cite{Ko1}, \cite{Ko2}, \cite{GST}).  Each of the six types of arcs (called tunnels) has the property that, once a regular neighborhood of the arc is added, then, up to isotopy, the regular neighborhood no longer depends on which $2$-bridge knot we started with -- indeed, we could have started with the unknot.   See Figure \ref{fig:twobridge}. In particular, there is a meridian of the complementary handlebody that intersects a meridian disk dual to the knot in a single point.  Following Lemma \ref{lemma:remstab} we then have

\begin{cor}  \label{cor:twobridge} Suppose $P$ is any bridge surface for a $2$-bridge knot $K \subset S^3$.  If $K$ is removable with respect to $P$, then $P$ is stabilized.
\end{cor}

In the proof of our main theorem we will use the following already known results. 

\begin{lemma} \cite[Lemma 3.1]{STo} \label{lem:bddcompressible}
 Let $A$ be a handlebody and let $(T, \bdd T) \subset (A, \bdd A)$ be a collection of bridges in $A$. Suppose $F$ is a properly embedded surface in $A$ transverse to $T$ that is not a union of unpunctured disks, once-punctured disks and twice-punctured spheres.  If $F_T$ is incompressible in $A_T$ then $\bdd F  \neq \emptyset$ and $F_T$ is $\bdd$-compressible.
 \end{lemma}

\begin{lemma} \cite[Lemma 3.6]{STo} \label{lem:disjointdisk}
Suppose $P$ and $Q$ are disjoint bridge surfaces for a link $K \subset M$, decomposing $M$ as $A \cup_P B$ and $X \cup_Q Y$ respectively.
Suppose furthermore that $Q_K\subset A_K$ and $P_K$ has a c-disk in $A_K$ that is disjoint from $Q_K$, 
    then either $P_K$ is c-strongly compressible or $M=S^3$ and $K$ 
    is empty or the unknot.  
 \end{lemma}

\begin{thm} \cite[Corollary 6.7]{STo} \label{thm:essential}
   Suppose $P$ and $Q$ are bridge surfaces for a link $K \subset M$ and $P_K$ and $Q_K$ are
   both c-weakly incompressible in $M_K$. If there is no incompressible Conway sphere for $K$ in $M$
then $P_K$ can be properly isotoped so that $P_K$ and $Q_K$ intersect in a non-empty collection of curves that are essential on both surfaces.
\end{thm}

\begin{thm} \cite{To} \label{thm:essentialexists}
  Suppose, for a link $K \subset M$, $M$ contains a c-strongly compressible bridge surface $Q$ that is not stabilized, meridionally stabilized or perturbed. Then either 
  \begin{itemize}
  \item $M$ contains a surface $F$ transverse to $K$ so that $F_K$ is essential in $M_K$ or
  \item $K$ contains a component $K_0$ that is removable. 
  \end{itemize}
  \end{thm}

\section{Unique Bridge Surface} \label{S:UniqueBridge}

Now we will focus our attention on two-bridge links in the 3-sphere. That is, for the rest of the paper, assume $S^3=A \cup_P B = X \cup_Q Y$, $K$ is in bridge position with respect to both $P$ and $Q$ and $P_K$ is a four times punctured sphere. In particular,  henceforth $A$ will be a ball that intersects $K$ in two trivial arcs. The ultimate goal is to show that if $K$ is non-trivial (i. e. neither the unknot nor the unlink of two components) and $Q_K$ is not  stabilized, meridionally stabilized or perturbed, then $Q_K$ is also a 4-times punctured sphere properly isotopic to $P_K$. We will 
  use the following technical lemma and its corollary.

\begin{lemma} \label{lem:parallel}
 Suppose $F_K$ is a connected splitting surface that is properly embedded in $A$, so $A_K = U_K \cup _{F_K} V_K$. Further assume $\bdd F$ consists of curves that are essential in $P_K$, $F_K$ is c-incompressible in $V_K$, but there is a $\bdd$-compressing disk for $F_K$ that lies in $V_K$.  Then $F_K$ is parallel to a subset of $P_K$ through $V_K$.  In particular $F_K$ is either an annulus or a twice punctured disk.
 \end{lemma}   

 \begin{proof}
     
Let $E\subset V_K$ be the $\bdd$-compressing disk for $F_K$. Let $\sigma = E \cap P_K$ and 
     note that $\sigma$ must be an essential arc on $P_K-F_K$ as 
     otherwise $F_K$ would be compressible in $V_K$. There are two cases to 
     consider.
     
     First suppose that both endpoints of $\sigma$ lie on the same component 
     of $\bdd F$; call this component $f$. As $f$ is an essential 
     curve on the 4-times punctured sphere $P_K$, it bounds two twice 
     punctured disks on $P_K$, let $P'$ be the twice 
     punctured disk containing $\sigma$. 
     A regular neighborhood of $P' \cup E$ consists of a copy 
     of $P'$ and two once punctured disks, $D'$ 
     and $D''$, whose boundaries lie on $F_K$. As $F_K$ is 
     c-incompressible in $V_K$, $D'$ and $D''$ each also bound once-punctured disks in $F_K$.  Moreover, these disks must be parallel to the once-punctured 
     disks on $F_K$, since twice-punctured spheres in a handlebody can only cut off trivial arcs from trivial arcs (cf \cite[Lemma 3.2]{STo}). Combining these parallelisms with the boundary 
     compression gives a parallelism betweet $F_K$ and $P'$.
     
     Suppose, on the other hand, that the two endpoints of $\sigma$ lie on different components of 
     $\bdd F$, say $f$ and $f'$. As $f$ and $f'$ are disjoint and essential in the 4-times punctured 
     sphere $P_K$,  $f$ and $f'$ must cobound an annulus $N$ on $P_K$ and $\sigma \subset 
     N$. A regular neighborhood of $N \cup E$ then consists of a copy 
     of $N$ and a disk $D$ whose boundary lies on $F_K$. As $F_K$ is 
     incompressible in $V_K$, $\bdd D$ also bounds a disk in $F_K$, a disk that is parallel to $D$ in $A_K$, since $A_K$ is irreducible. Combining  this parallelism with the boundary compression gives the desired  parallelism between $F_K$ and $N$. \end{proof}

\begin{cor} \label{cor:parallel}
    Suppose $F_K$ is a c-incompressible connected splitting surface, not an unpunctured disk, that is properly embedded in $A$, and suppose $\bdd F$ consists of curves that are essential in $P_K$.  Then $F_K$ is $P_K$-parallel.
 \end{cor}   
 
 \begin{proof}
 $F_K$ can't be a once-punctured disk, since its boundary also bounds a twice-punctured disk in $P_K$.    Since it's c-incompressible, it's incompressible, so by Lemma \ref{lem:bddcompressible}, $F_K$ must be boundary-compressible. The result follows by Lemma \ref{lem:parallel}
 \end{proof}

\begin{thm} \label{thm:cweaklyincomp}
    Let $K \subset S^3$ be a two bridge link (not a trivial knot or link) with respect to a bridge surface $P \cong S^2 \subset S^3$.  Any c-weakly incompressible bridge surface for $(S^3,K)$ is properly isotopic to $P_K$.

\end{thm}
\begin{proof}
Suppose $Q$ is a c-weakly incompressible bridge surface, so $S^3 = A \cup_P B = X \cup_Q Y$.  $P$ is also c-weakly incompressible.  Indeed, disjoint essential curves in the $4$-punctured sphere $P$ are necessarily parallel in $P_K$, and so a c-strong compressing pair would provide a splitting sphere for $K$, contradicting the assumption that $K$ is not a trivial link.  By Theorem \ref{thm:essential} we may isotope $P_K$ so that $P_K \cap Q_K \neq \emptyset$ and all curves of  $P_K \cap Q_K$ are essential on both $P_K$ and $Q_K$. Furthermore assume that the number of components of intersection $|P_K \cap Q_K|$ is minimal under these restrictions. We will denote by $Q^A_K$ and $Q^B_K$ the surfaces $Q_K \cap A$ and $Q_K \cap B$ respectively. Similarly we will denote by $P^X_K$  and $P^Y_K$ the surfaces $P_K \cap X$ and $P_K \cap Y$.
\medskip

\textbf{Claim 1:  At least one of $Q_K^A$ or $Q_K^B$ has a $P$-parallel component.}
\medskip

$Q_K$ is not a twice-punctured sphere, since $K$ is not the unknot.  Thus there are c-disks for $Q_K$ in both $X$ and $Y$.  First we will reduce to the case that there are c-disks for $Q_K$ in both $X$ and $Y$ that are both disjoint from $P_K$. 

If there aren't such c-disks, then, with no loss of generality, there is a c-disk $D^*_Y \subset Y$ for $Q_K$ so that $|P \cap D^*_Y| > 0$ is  minimal among all c-disks for $Q_K$ in $Y$.  If the intersection contains any simple closed curves, let $\alpha$ be an innermost one on $D^*_Y$ bounding a possibly punctured
disk $D^*_{\alpha} \subset D^*_Y$. If $\alpha$ were inessential in $P$, then a c-disk with fewer intersection curves could have been found, so $\alpha$ is essential in $P_K$.  Note that as $P_K$ is a 4-times punctured sphere and all curves of $P_K \cap Q_K$ are essential in $P_K$, all the curves must be parallel on $P_K$ and are all also parallel to $\alpha$. Let $N \subset P_K$ be the annulus
between $\alpha$ and an adjacent curve of $P_K \cap Q_K$. Then by 
slightly isotoping the possibly punctured disk $N \cup D^*_{\alpha}$ we obtain c-disk for $Q_K$ that is disjoint from $P_K$ contradicting the choice of $D^*_Y$. Thus we may assume that $D^*_Y 
\cap P_K$ consists only of arcs. An arc of $D^*_Y  \cap P_K$ that is outermost on $D^*_Y$ cuts off a disk in $Y_K$ that $\bdd$-compresses $Q_K^A$, say, to $P_K$. By Lemma \ref{lem:parallel}, $Q_K^A$ has a component that is $P_K$-parallel, establishing the claim in this case.

So now assume that there are c-disks $D^*_Y \subset Y$ and $D^*_X \subset X$ for $Q_K$ and both are disjoint from $P_K$.  If one disk lies in $A$ and the other in $B$, then the disks would have disjoint boundaries, contradicting the assumption that $Q_K$ is c-weakly incompressible. So these c-disks must lie on the same side of $P_K$.  Suppose without loss that they both lie in $A$.  Then, since $Q_K$ is c-weakly incompressible, $Q_K^B$ must be c-incompressible in $B$. But by Corollary \ref{cor:parallel}, this implies that $Q_K^B$ has a $P_K$-parallel component, again establishing the claim.

\medskip

Following the claim, suppose with no loss of generality that $Q_K^A$ has a $P_K$-parallel component. In this case $Q_K^A$  must be connected, for otherwise a  component of $Q_K^A$ could be isotoped across $P_K$ reducing $|P_K \cap Q_K|$. As all components of $P_K-Q_K$ are annuli or twice punctured disks, $Q_K^A$ is also either an annulus or a twice punctured disk. Without loss of
generality, assume $Q_K^A$ is parallel to $P_K^X$ (through the region $A \cap X$).  See Figure \ref{fig:Venn}.

 \begin{figure}[tbh]
    \centering
    \includegraphics[scale=0.6]{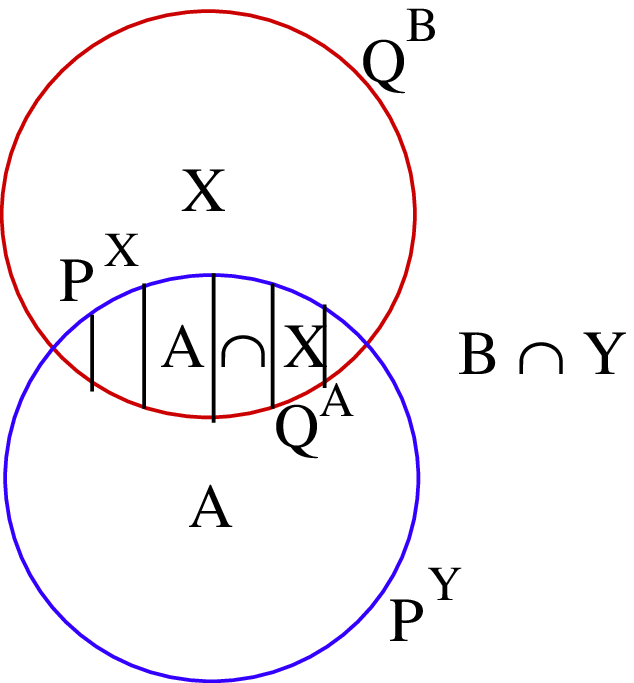}
    \caption{} \label{fig:Venn}
    \end{figure}

Suppose $Q_K^B$ were c-compressible into $Y_K$ with a c-disk 
$D^*$. Isotope $P_K^X$ across $Q_K^A$ so that $P_K \subset Y_K-D^*$. By Lemma 
\ref{lem:disjointdisk} this would imply that $Q_K$ is c-strongly compressible, a contradiction to our hypothesis. Thus we conclude that $Q_K^B$ is either c-incompressible or c-compresses only into $X_K$. A similar argument  with the roles of $P$ and $Q$ switched shows that $P_K^Y$ does not c-compress into $B_K$.

\medskip

\textbf{Case 1:  $Q_K^A$ consists of a single $P$-parallel twice punctured disk.}

This in particular implies that both $P_K^X$ and $P_K^Y$ consist of single twice-punctured disks.

Suppose first that $P_K^Y$ is c-compressible in $Y$  with 
c-disk $D^*$. As already shown $D^*\subset A$.  Without loss of generality we may assume $\bdd D^* = (P_K \cap 
Q_K)$ so $D^*$ is also a c-disk for $Q_K^A$ 
lying in $Y_K$. As $Q_K$ is c-weakly incompressible, $Q_K^B$ is either 
c-incompressible or also c-compresses in $Y_K$. As we have already 
eliminated the later option, $Q_K^B$ must be c-incompressible and so by Corollary 
\ref{cor:parallel} $Q_K^B$ is $P_K$-parallel. Thus
$Q_K^B$ is a twice punctured disk so $Q_K$ is 
also a 4-times punctured sphere. 
In summary, if $P_K^Y$ is c-compressible, then $Q_K$ is also a 
4-times punctured sphere and $Q_K^B$ is c-incompressible.  So, by 
possibly switching the names of $P$ and $Q$, we may henceforth assume that $P_K^Y$ is 
c-incompressible.

As $P_K^Y$ is c-incompressible, by Lemma \ref{lem:bddcompressible} it must 
be $\bdd$-compressible. Let $E$ be the boundary compressing disk and 
note that $E \cap Q_K$ is an arc essential on $Q_K-P_K$ as otherwise 
$P_K^Y$ would be compressible. Thus, by changing our point of view, we 
can conclude that $Q_K^A$ or $Q_K^B$ is $\bdd$-compressible in $A$ or $B$ respectively to $P_K^Y$.

Suppose $Q_K^B$ is $\bdd$-compressible to $P_K^Y$.  As $Q_K^B$ 
is c-incompressible in $Y$, Lemma \ref{lem:parallel} 
implies that $Q_K^B$ is parallel to $P_K^Y$. Combining this with 
parallelism between $Q_K^A$ and $P_K^X$ gives the desired isotopy between $P_K$ 
and $Q_K$.

Suppose $Q_K^A$ is $\bdd$-compressible into $P_K^Y$. 
Since $P_K^Y$ is a c-incompressible splitting surface for $Y$, it follows from Lemma \ref{lem:parallel} that $P_K^Y$ is parallel to $Q_K^A$, i.e.  $Q_K^A$ is 
 isotopic to both $P_K^X$ and $P_K^Y$.  In particular $P_K$ can be properly
 isotoped to lie in either $X_K$ or $Y_K$. By Lemma \ref {lem:disjointdisk}
 this implies that $Q_K^B$ must be 
 c-incompressible in $B_K$, for if $Q_K^B$ has a c-disk lying in $X_K$ (say) 
 we could isotope $P_K$ to lie in $X_K$ and be disjoint from this 
 c-disk. By Corollary \ref{cor:parallel} this implies that 
 $Q_K^B$ is parallel to one of $P_K^X$ or $P_K^Y$. As $Q_K^A$ is 
 parallel to both $P_K^X$ and $P_K^Y$ we conclude that $P_K$ and 
 $Q_K$ are properly isotopic.

\medskip

\textbf{Case 2:  $Q_K^A$ is a single $P_K^X$-parallel annulus.}

We will show, by contradicition, that this case does not arise.  In this situation $P_K^X$ is a single
annulus and $P_K^Y$ consists of two twice-punctured disks. See Figure \ref{fig:parannuli}. Recall 
that we have already shown that $Q_K^B$ is c-incompressible in $Y_K$ 
and $P_K^Y$ is c-incompressible in $B_K$.

 \begin{figure}[tbh]
    \centering
    \includegraphics[scale=0.6]{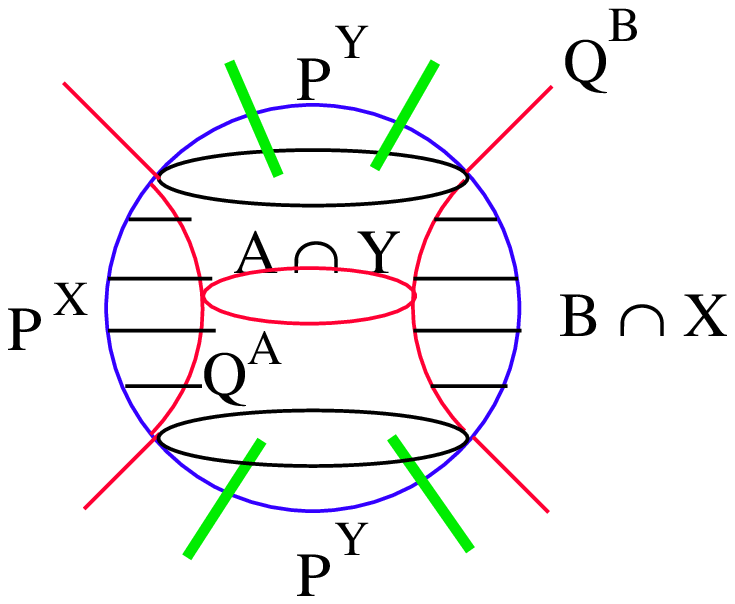}
    \caption{} \label{fig:parannuli}
    \end{figure}

Suppose (towards a contradiction) that $P_K^Y$ is c-incompressible in $Y$. By 
Lemma \ref{lem:bddcompressible} it must be boundary compressible. As in the previous 
case the $\bdd$-compressing 
disk is incident to $Q_K-P_K$ in an essential arc, i.e. one of $Q_K^A$ or 
$Q_K^B$ is $\bdd$-compressible to $P_K^Y$. The annulus $Q_K^A$ can't $\bdd$-compresses 
to $P_K^Y$, since its boundary components are on different components of $P_K^Y$.  On the other hand, if a component of $Q_K^B$ $\bdd$-compresses to $P_K^Y$, by Corollary \ref{cor:parallel} it 
follows  that $Q_K^B$ has a twice punctured disk component parallel to one of 
the two components of $P_K^Y$. In this case $|P_K \cap Q_K|$ can be decreased by 1, and this contradicts the  minimality assumption. We conclude that $P_K^Y$ must be c-compressible 
in the complement of $Q_K$.  

Suppose $D^*$ is a c-disk for $P_K^Y$ in the complement of  $Q_K$.  Necessarily $D^* \subset A_K$ and we may as well take $\bdd D^*$ to be one of the circles $P \cap Q$. Then $D^*$ is also a c-disk for $Q_K^A$  lying in $A \cap Y$. By c-weak incompressibility of $Q_K$ any c-disks for $Q_K^B$ 
would have to lie in $Y_K$.  But we established in the beginning of this case that
this is not possible so $Q_K^B$ is in fact c-incompressible. By 
Lemma \ref{lem:bddcompressible} $Q_K^B$ must be boundary compressible.
As we already saw, if the boundary compression is to $P_K^Y$, the 
intersection $P_K \cap Q_K$ can be reduced, so $Q_K^B$ must be boundary 
compressible to the annulus $P_K^X$.  It follows then from Lemma \ref{lem:parallel} that $Q_K^B$, like $Q_K^A$, is an annulus parallel to $P_K^X$.  Then $Q_K$ is a torus that is disjoint from $K$ and so it cannot be a bridge surface, a contradiction.  
\end{proof}

\begin{cor}  Suppose $K$ is a knot in $S^3$, $2$-bridge with respect to the bridge surface $P \cong S^2$, and $K$ is not the unknot.  Suppose $Q$ is any other bridge surface for $K$.  Then either
\begin{itemize}
\item $Q$ is stabilized
\item $Q$ is meridionally stabilized
\item $Q$ is perturbed or
\item $Q$ is properly isotopic to $P$.
\end{itemize}
\end{cor}

\begin{proof}  If $Q$ is c-weakly incompressible then Theorem \ref{thm:cweaklyincomp} shows that $Q$ is properly isotopic to $P$. If $Q$ is c-strongly compressible, Theorem \ref{thm:essentialexists} says that either $Q$ is stabilized, meridionally stabilized or perturbed, or $K$ is removable with respect to the bridge surface $Q$, or there is a surface $F$ transverse to $K$ so that $F_K$ is essential in $S^3_K$.  The last possibility does not occur for $2$-bridge knots (see \cite{HT}).  Corollary \ref{cor:twobridge} shows that if $K$ is removable with respect to $Q$, then $Q$ is stabilized.
\end{proof}

 \bibliography{mybibliounique}
 \bibliographystyle{plain}

\end{document}